\def\author{S. Kleiman and A. Thorup}
\def\title{Conormal Geometry of Maximal Minors}
\def\date{12 May 1998}
\def\TheMagstep{\magstep1}      
\def\PaperSize{letter}          
\def\hdfnt{\smc}
\let\:=\colon \let\ox=\otimes \let\x=\times 
\let\wt=\widetilde \let\?=\overline \let\into=\hookrightarrow
\let\To=\longrightarrow 
\def\IV{{\bf V}} \def\IP{{\bf P}} \def\IC{{\bf C}}

\def\mcl#1{\expandafter\def\csname c#1\endcsname{{\cal#1}}}
 \mcl E \mcl F \mcl I \mcl M \mcl N \mcl O
\def\mop#1 {\expandafter\def\csname#1\endcsname
  {\mathop{\rm#1}\nolimits}}
\mop Proj \mop Spec \mop Supp \mop Bl \mop Image \mop Sing \mop Projan
\def\hgt{\mathop{\rm ht}\nolimits}  \def\as{\mathop{\ell}\nolimits}
\def\trdg{\mathop{\rm tr.deg}\nolimits}
\def\Ext{\mathop{{\cal E}\!{\it xt}}\nolimits}
\def\Hom{\mathop{{\cal H}\!{\it om}}\nolimits}
\def\and{\hbox{\rm\ and }} \def\where{\hbox{\rm\ where }}

\parskip=0pt plus 1.75pt \parindent10pt
\hsize28pc 
\vsize43pc 
\abovedisplayskip 4pt plus3pt minus1pt
\belowdisplayskip=\abovedisplayskip
\abovedisplayshortskip 2.5pt plus2pt minus1pt
\belowdisplayshortskip=\abovedisplayskip

\def\today{\ifcase\month\or     
 January\or February\or March\or April\or May\or June\or
 July\or August\or September\or October\or November\or December\fi
 \space\number\day, \number\year}
\nopagenumbers
\headline={\hdfnt
  \ifnum\pageno=1\firstheadline
  \else
    \ifodd\pageno\oddheadline
    \else\evenheadline\fi
  \fi
}
\expandafter\ifx\csname date\endcsname\relax \let\dato=\today
	    \else\let\dato=\date\fi
\let\firstheadline\hfill
\def\oddheadline{\rlap{\dato}\hfil\hdtitle\hfil\llap{\folio}}
\def\evenheadline{\rlap{\folio} \hfil\hdauthor\hfil\llap{\dato}}

\def\TRUE{TRUE} 
\ifx\DoublepageOutput\TRUE \def\TheMagstep{\magstep0} \fi
\mag=\TheMagstep

\newskip\vadjustskip \vadjustskip=0.5\normalbaselineskip
\def\centertext
 {\hoffset=\pgwidth \advance\hoffset-\hsize
  \advance\hoffset-2truein \divide\hoffset by 2\relax
  \voffset=\pgheight \advance\voffset-\vsize
  \advance\voffset-2truein \divide\voffset by 2\relax
  \advance\voffset\vadjustskip
 }
\newdimen\pgwidth\newdimen\pgheight
\def\letter{letter}\def\AFour{AFour}
\ifx\PaperSize\letter
 \pgwidth=8.5truein \pgheight=11truein
 \message{- Got a paper size of letter.  }\centertext
\fi
\ifx\PaperSize\AFour
 \pgwidth=210truemm \pgheight=297truemm
 \message{- Got a paper size of AFour.  }\centertext
\fi

 \font\twelvebf=cmbx12          
 \font\smc=cmcsc10              
\catcode`\@=11          
\def\eightpoint{\eightpointfonts
 \setbox\strutbox\hbox{\vrule height7\p@ depth2\p@ width\z@}%
 \eightpointparameters\eightpointfamilies
 \normalbaselines\rm
 }
\def\eightpointparameters{%
 \normalbaselineskip9\p@
 \abovedisplayskip9\p@ plus2.4\p@ minus6.2\p@
 \belowdisplayskip9\p@ plus2.4\p@ minus6.2\p@
 \abovedisplayshortskip\z@ plus2.4\p@
 \belowdisplayshortskip5.6\p@ plus2.4\p@ minus3.2\p@
 }
\newfam\smcfam
\def\eightpointfonts{%
 \font\eightrm=cmr8 \font\sixrm=cmr6
 \font\eightbf=cmbx8 \font\sixbf=cmbx6
 \font\eightit=cmti8
 \font\eightsmc=cmcsc8
 \font\eighti=cmmi8 \font\sixi=cmmi6
 \font\eightsy=cmsy8 \font\sixsy=cmsy6
 \font\eightsl=cmsl8 \font\eighttt=cmtt8}
\def\eightpointfamilies{%
 \textfont\z@\eightrm \scriptfont\z@\sixrm  \scriptscriptfont\z@\fiverm
 \textfont\@ne\eighti \scriptfont\@ne\sixi  \scriptscriptfont\@ne\fivei
 \textfont\tw@\eightsy \scriptfont\tw@\sixsy \scriptscriptfont\tw@\fivesy
 \textfont\thr@@\tenex \scriptfont\thr@@\tenex\scriptscriptfont\thr@@\tenex
 \textfont\itfam\eightit        \def\it{\fam\itfam\eightit}%
 \textfont\slfam\eightsl        \def\sl{\fam\slfam\eightsl}%
 \textfont\ttfam\eighttt        \def\tt{\fam\ttfam\eighttt}%
 \textfont\smcfam\eightsmc      \def\smc{\fam\smcfam\eightsmc}%
 \textfont\bffam\eightbf \scriptfont\bffam\sixbf
   \scriptscriptfont\bffam\fivebf       \def\bf{\fam\bffam\eightbf}%
 \def\rm{\fam0\eightrm}%
 }
\def\vfootnote#1{\insert\footins\bgroup
 \eightpoint 
 \interlinepenalty\interfootnotelinepenalty
  \splittopskip\ht\strutbox 
  \splitmaxdepth\dp\strutbox \floatingpenalty\@MM
  \leftskip\z@skip \rightskip\z@skip \spaceskip\z@skip \xspaceskip\z@skip
  \textindent{#1}\footstrut\futurelet\next\fo@t}

 \newdimen\fullhsize \newbox\leftcolumn
 \def\fulline{\hbox to \fullhsize}
\def\doublepageoutput
{\let\lr=L
 \output={\if L\lr
          \global\setbox\leftcolumn=\columnbox \global\let\lr=R%
        \else \doubleformat \global\let\lr=L\fi
        \ifnum\outputpenalty>-20000 \else\dosupereject\fi}%
 \def\doubleformat{\shipout\vbox{%
     \ifx\PaperSize\AFour
           \fulline{\hfil\box\leftcolumn\hfil\columnbox\hfil}%
     \else
           \fulline{\hfil\hfil\box\leftcolumn\hfil\columnbox\hfil\hfil}%
     \fi             }%
}
 \def\columnbox{\vbox
   {\ifvoid255\headline={\hfil}\nopagenumbers\else
    \makeheadline\pagebody\makefootline\advancepageno\fi}%
   }%
 \fullhsize=\pgheight \hoffset=-1truein
 \voffset=\pgwidth \advance\voffset-\vsize
  \advance\voffset-2truein \divide\voffset by 2
  \advance\voffset\vadjustskip
 \def\firstheadline{\leftline{alg-geom/9708018}}
 
\ifx\FirstPageOnRight\TRUE 
 \null\vfill\nopagenumbers\eject\pageno=1\relax
\fi
}
\ifx\DoublepageOutput\TRUE \doublepageoutput \fi


\newcount\sctno
\def\sctn#1\par
  {\removelastskip\vskip0pt plus3\normalbaselineskip
 \penalty-250
  \vskip0pt plus-3\normalbaselineskip \bigskip\smallskip
  \centerline{#1}\medskip
  }
\def\sct#1.{\sctno=#1\relax\sctn#1.}

\def\Cs#1){\(\number\sctno.#1)}
\def\part#1 {\par\(#1)\enspace\ignorespaces}
\def\qis{\quad\ignorespaces}
\def\qbi{\quad\begingroup\it}
\def\dsc#1 #2.{\medbreak\Cs#1) {\it #2.}\qis}
\def\endpro{\endgroup\medbreak}
\def\lem#1 {\medbreak\Cs#1) {\smc Lemma.}\qbi}
\def\pro#1 {\medbreak{\smc#1.}\qbi}

 \newcount\refno \refno=0        \def\NoKey{*!*}
 \def\MakeKey{\advance\refno by 1 \expandafter\xdef
  \csname\TheKey\endcsname{{\number\refno}}\NextKey}
 \def\NextKey#1 {\def\TheKey{#1}\ifx\TheKey\NoKey\let\next\relax
  \else\let\next\MakeKey \fi \next}
 \def\RefKeys #1\endRefKeys{\expandafter\NextKey #1 *!* }
 \def\SetRef#1 #2,{\hang\llap
  {\csname#1\endcsname.\enspace}{\smc #2},}
 \newbox\keybox \setbox\keybox=\hbox{8.\enspace}
 \newdimen\keyindent \keyindent=\wd\keybox
\def\references{
  \bgroup   \frenchspacing   \eightpoint
   \parindent=\keyindent  \parskip=\smallskipamount
   \everypar={\SetRef}\par}
\def\endreferences{\egroup}

 \def\serial#1#2{\expandafter\def\csname#1\endcsname ##1 ##2
        ##3.{\unskip\ {\it #2\/} {\bf##1} (##2), ##3.}} 
 \serial{tams}{Trans. Amer. Math. Soc.}
 \serial{comp}{Comp. Math.}
 \serial{conm}{Contemp. Math.}
 \serial{ja}{J. Algebra}
 \serial{slnm}{Springer Lecture Notes in Math}

\def\UThin{\penalty\@M \thinspace\ignorespaces}
\def\(#1){{\let~=\UThin\rm(#1)}}
\def\relaxnext@{\let\next\relax}
\def\cite#1{\relaxnext@
 \def\nextiii@##1,##2\end@{\unskip\space{\rm[\SetKey{##1},\let~=\UThin##2]}}%
 \in@,{#1}\ifin@\def\next{\nextiii@#1\end@}\else
 \def\next{{\rm[\SetKey{#1}]}}\fi\next}
\newif\ifin@
\def\in@#1#2{\def\in@@##1#1##2##3\in@@
 {\ifx\in@##2\in@false\else\in@true\fi}%
 \in@@#2#1\in@\in@@}
\def\SetKey#1{{\bf\csname#1\endcsname}}

\catcode`\@=12  


\def\hdauthor{kleiman and thorup}
\def\hdtitle{conormal geometry of maximal minors}
\def\oddheadline{\hfil\hdtitle\hfil\llap{\folio}}
\def\evenheadline{\rlap{\folio}\hfil\hdauthor\hfil\llap{\dato}}
\RefKeys
  B BH BR BV GM GK KR K76 K85 K94 K98 KL KT Lo84 M98 Mat Mat90 R89 T
 \endRefKeys

\let\at=@\catcode`\@=11

 \def\activeat#1{\csname @#1\endcsname}
 \def\def@#1{\expandafter\def\csname @#1\endcsname}
 {\catcode`\@=\active \gdef@{\activeat}}

\let\ssize\scriptstyle
\newdimen\ex@	\ex@.2326ex

 \def\requalfill{\cleaders\hbox{$\mkern-2mu\mathord=\mkern-2mu$}\hfill
  \mkern-6mu\mathord=$}
 \def\eqfill{$\m@th\mathord=\mkern-6mu\requalfill}
 \def\deffill{\hbox{$:=$}$\m@th\mkern-6mu\requalfill}
 \def\fiberbox{\hbox{$\vcenter{\hrule\hbox{\vrule\kern1ex
     \vbox{\kern1.2ex}\vrule}\hrule}$}}

 \font\arrfont=line10
 \def\Swarrow{\vcenter{\hbox{$\swarrow$\kern-.26ex
    \raise1.5ex\hbox{\arrfont\char'000}}}}

 \newdimen\arrwd
 \newdimen\minCDarrwd \minCDarrwd=2.5pc
 	
 \def\findarrwd#1#2#3{\arrwd=#3%
  \setbox\z@\hbox{$\ssize\;{#1}\;\;$}%
  \setbox\@ne\hbox{$\ssize\;{#2}\;\;$}%
  \ifdim\wd\z@>\arrwd \arrwd=\wd\z@\fi
  \ifdim\wd\@ne>\arrwd \arrwd=\wd\@ne\fi}
 \newdimen\arrowsp\arrowsp=0.375em
 \def\findCDarrwd#1#2{\findarrwd{#1}{#2}{\minCDarrwd}
    \advance\arrwd by 2\arrowsp}
 \newdimen\minarrwd 
 \setbox\z@\hbox{$\longrightarrow$} \minarrwd=\wd\z@

 \def\harrow#1#2#3#4{{\minarrwd=#1\minarrwd%
   \findarrwd{#2}{#3}{\minarrwd}\kern\arrowsp
    \mathrel{\mathop{\hbox to\arrwd{#4}}\limits^{#2}_{#3}}\kern\arrowsp}}
 \def@]#1>#2>#3>{\harrow{#1}{#2}{#3}\rightarrowfill}
 \def@>#1>#2>{\harrow1{#1}{#2}\rightarrowfill}
 \def@<#1<#2<{\harrow1{#1}{#2}\leftarrowfill}
 \def@={\harrow1{}{}\eqfill}
 \def@:#1={\harrow1{}{}\deffill}
 \def@ N#1N#2N{\vCDarrow{#1}{#2}\UpDownarrow}
 \def\UpDownarrow{\uparrow\,\Big\downarrow}

 \def@.{\ifodd\row\relax\harrow1{}{}\hfill
   \else\vCDarrow{}{}.\fi}
 \def@|{\vCDarrow{}{}\Vert}
 \def@ V#1V#2V{\vCDarrow{#1}{#2}\downarrow}
 \def@ A#1A#2A{\vCDarrow{#1}{#2}\uparrow}
 \def@(#1){\arrwd=\csname col\the\col\endcsname\relax
   \hbox to 0pt{\hbox to \arrwd{\hss$\vcenter{\hbox{$#1$}}$\hss}\hss}}

 \def\squash#1{\setbox\z@=\hbox{$#1$}\finsm@@sh}
\def\finsm@@sh{\ifnum\row>1\ht\z@\z@\fi \dp\z@\z@ \box\z@}

 \newcount\row \newcount\col \newcount\numcol \newcount\arrspan
 \newdimen\vrtxhalfwd  \newbox\tempbox

 \def\innernewdimen{\alloc@1\dimen\dimendef\insc@unt}
 \def\measureinit{\col=1\vrtxhalfwd=0pt\arrspan=1\arrwd=0pt
   \setbox\tempbox=\hbox\bgroup$}
 \def\setinit{\col=1\hbox\bgroup$\ifodd\row
   \kern\csname col1\endcsname
   \kern-\csname row\the\row col1\endcsname\fi}
 \def\findvrtxhalfsum{$\egroup
  \expandafter\innernewdimen\csname row\the\row col\the\col\endcsname
  \global\csname row\the\row col\the\col\endcsname=\vrtxhalfwd
  \vrtxhalfwd=0.5\wd\tempbox
  \global\advance\csname row\the\row col\the\col\endcsname by \vrtxhalfwd
  \advance\arrwd by \csname row\the\row col\the\col\endcsname
  \divide\arrwd by \arrspan
  \loop\ifnum\col>\numcol \numcol=\col%
     \expandafter\innernewdimen \csname col\the\col\endcsname
     \global\csname col\the\col\endcsname=\arrwd
   \else \ifdim\arrwd >\csname col\the\col\endcsname
      \global\csname col\the\col\endcsname=\arrwd\fi\fi
   \advance\arrspan by -1 %
   \ifnum\arrspan>0 \repeat}
 \def\setCDarrow#1#2#3#4{\advance\col by 1 \arrspan=#1
    \arrwd= -\csname row\the\row col\the\col\endcsname\relax
    \loop\advance\arrwd by \csname col\the\col\endcsname
     \ifnum\arrspan>1 \advance\col by 1 \advance\arrspan by -1%
     \repeat
    \squash{\mathop{
     \hbox to\arrwd{\kern\arrowsp#4\kern\arrowsp}}\limits^{#2}_{#3}}}
 \def\measureCDarrow#1#2#3#4{\findvrtxhalfsum\advance\col by 1%
   \arrspan=#1\findCDarrwd{#2}{#3}%
    \setbox\tempbox=\hbox\bgroup$}
 \def\vCDarrow#1#2#3{\kern\csname col\the\col\endcsname
    \hbox to 0pt{\hss$\vcenter{\llap{$\ssize#1$}}%
     \Big#3\vcenter{\rlap{$\ssize#2$}}$\hss}\advance\col by 1}

 \def\setCD{\def\harrow{\setCDarrow}%
  \def\\{$\egroup\advance\row by 1\setinit}
  \m@th\lineskip3\ex@\lineskiplimit3\ex@ \row=1\setinit}
 \def\endsetCD{$\egroup}
 \def\measure{\bgroup
  \def\harrow{\measureCDarrow}%
  \def\\##1\\{\findvrtxhalfsum\advance\row by 2 \measureinit}%
  \row=1\numcol=0\measureinit}
 \def\endmeasure{\findvrtxhalfsum\egroup}

\newbox\CDbox \newdimen\sdim

 \newcount\savedcount
 \def\CD#1\endCD{\savedcount=\count11%
   \measure#1\endmeasure
   \vcenter{\setCD#1\endsetCD}%
   \global\count11=\savedcount}

 \catcode`\@=\active
{\leftskip=0pt plus1fill \rightskip=\leftskip
 \obeylines
 \leavevmode \medskip
 {\twelvebf \title
 } \medskip
 \footnote{}{\noindent %
 Subj-class: 13C15 (Primary) 14N05, 14B05 (Secondary).}
 Steven Kleiman\footnote{$^{1}$}{%
    Supported in part by NSF grant 9400918-DMS.}
 {\eightpoint\it\medskip
 Deptartment of Mathematics, Room {\sl 2-278} MIT,
 {\sl77} Mass Ave, Cambridge, MA {\sl02139-4307}, USA
 \rm E-mail: \tt kleiman\at  math.mit.edu \medskip
 } and \medskip
 \rm Anders Thorup\footnote{$^{2}$}{%
 Supported in part by the Danish Natural Science Research Council, grant %
 11-7428.}
 {\eightpoint\it\medskip
 Matematisk Afdeling, K\o benhavns Universitet,
 Universitetsparken {\sl5}, DK-{\sl2100} K\o benhavn \O, Danmark
 \rm E-mail: \tt thorup\at math.ku.dk \medskip\smallskip
 \rm \dato \bigskip
 }
}
{\parindent=1.5\parindent \narrower \eightpoint
 Let $A$ be a Noetherian local domain, $N$ be a finitely generated
torsion-free module, and $M$ a proper submodule that is generically
equal to $N$.  Let $A[N]$ be an arbitrary graded overdomain of $A$
generated as an $A$-algebra by $N$ placed in degree 1.  Let $A[M]$ be
the subalgebra generated by $M$.  Set $C:=\Proj(A[M])$ and $r:=\dim C$.
Form the (closed) subset $W$ of $\Spec(A)$ of primes $\bf p$ where
$A[N]_{\bf p}$ is not a finitely generated module over $A[M]_{\bf p}$,
and denote the preimage of $W$ in $C$ by $E$.  We prove this: {\it
\(1)~$\dim E=r-1$ if either \(a)~$N$ is free and $A[N]$ is the symmetric
algebra, or \(b)~$W$ is nonempty and $A$ is universally catenary, and
\(2)~$E$ is equidimensional if \(a) holds and $A$ is universally
catenary.}

Our proof was inspired by some recent work of Gaffney and Massey, which
we sketch; they proved (2) when $A$ is the ring of germs of a
complex-analytic variety, and applied it to perfect a characterization
of Thom's A$_f$-condition in equisingularity theory.  From (1), we
recover, with new proofs, the usual height inequality for maximal minors
and an extension of it obtained by the authors in 1992.  From the
latter, we recover the authors' generalization to modules of B\"oger's
criterion for integral dependence of ideals.  Finally, we introduce an
application of (1), being made by the second author, to the geometry of
the dual variety of a projective variety, and use it to obtain an
interesting example where the conclusion of (1) fails and $A[N]$ is a
finitely generated module over $A[M]$.
  \par } 

\sct1. The Theorem and Applications

\dsc1 Introduction.  Our main result is the following theorem.  Its
proof occupies nearly all of the next section.

Let $A$ be a Noetherian local domain, $N$ be a finitely generated
torsion-free module, and $M$ a nonzero proper submodule.  Set
$X:=\Spec(A)$ and $Y:=\Supp(N/M)$.  Let $A[N]$ be an arbitrary graded
domain containing $A$ and generated as an $A$-algebra by $N$ placed in
degree 1.  Thus $A[N]$ either is the {\it Rees algebra} (that is, the
quotient of the symmetric algebra by its torsion) or is a quotient of
the Rees algebra by a homogeneous prime ideal that intersects $N$ in
$0$.  Let $A[M]$ be the subalgebra generated by $M$.  Set
$P:=\Proj(A[N])$, let $p\:P\to X$ denote the structure map, and set
$Z:=\IV(M\cdot A[N])$; so $Z$ is the subscheme of $P$ whose homogeneous
ideal is generated by $M$.  Set $C:=\Proj(A[M])$, let $c\:C\to X$ denote
the structure map, and set $E:=c^{-1}p(Z)$ and $F:=c^{-1}Y$.  Finally,
set $r:=\dim P$.

Note that $p(Z)\subset Y$ since, off $Y$, the ideal $M\cdot A[N]$ is
irrelevant; so $E\subset F$.  Note that $P$ and $C$ are integral, and
that $p\:P\to X$ and $c\:C\to X$ are surjective, being proper and being
dominating as $A\subset A[M]\subset A[N]$.  If $Y\neq X$, then
generically $M$ and $N$ are equal; whence, by (3.4)(ii) of \cite{KT},
        $$\dim C=r.\eqno\Cs1.1)$$

\pro Theorem Preserve the notation above, and assume $Y\neq X$.
 \part1 If either \(a)~$N$ is free, and $A[N]$ is the symmetric algebra,
\(b)~$A[N]$ is not a finitely generated $A[M]$-module, and $A$ is
universally catenary, or \(c)~$Z=p^{-1}p(Z)$ as sets, and $Z$ is
nonempty, or \(d)~$P$ has dimension $r$ at some point of $Z$, then
  $$\dim E=r-1 \and \dim F=r-1.$$
 Furthermore, if \(a) holds, then $p(Z)=Y$ and $E=F$.
 \part2 Assume either that $N$ is free, and $A[N]$ is the symmetric
algebra, or that $Z=p^{-1}p(Z)$ as sets.  If $A$ is universally
catenary, then $C$, $P$, and $E$ are biequidimensional.  \endpro

The theorem has applications in algebra and in geometry, which will be
discussed in this section.  In short, Part~(1) with Hypothesis~(a)
implies the usual height inequality for maximal minors, 
because $Y$ is defined by the zeroth Fitting ideal $I$ of $N/M$.  The
height inequality was given one of its first proofs by Buchsbaum and Rim
\cite{BR, 3.5} as an application of their theory of multiplicities of
submodules of free modules.  (See \cite{BV, Ch.~2} for a discussion of
other proofs.)  Following in their footsteps, but assuming that $p(Z)$
is the closed point, the authors recovered (1) with (a) and proved (1)
with (b)--(d) in (10.2) and (10.3) of \cite{KT}.  These four results are
recovered in the present article via new proofs; moreover, these
proofs are substantially shorter, simpler, and more direct than the old.
If $A$ is universally catenary, then (1) can be reduced to the case
where $p(Z)$ is the closed point by localizing at a generic point of
$p(Z)$; if $A$ is not universally catenary, then (1) appears to be new.

Part~(1) with Hypothesis~(b) provides a criterion for $N$ to be
``integrally dependent'' on $M$.  As such, (1) with (b) is the main
ingredient in the authors' generalization to modules \cite{KT, (10.9)}
of B\"oger's criterion, which in turn generalized to ideals not of
finite colength Rees's celebrated characterization of integral
dependence by multiplicity.  Part~(1) with (b) is therefore a main
ingredient in the work \cite{GK} of Gaffney and the first author, which
generalizes to modules Teissier's principle of specialization of
integral dependence, and applies it in equisingularity theory.

Part~(2) asserts the pure (graded) codimensionality of the extension
$I\cdot A[M]$ of the Fitting ideal $I$ to the Rees algebra $A[M]$,
without any assumption on the codimensionality of $I$ itself.
 (Graded (co)dimension is defined via chains of homogeneous primes, but
is equal to the usual notion, defined via chains of
arbitrary primes by Theorem~1.5.8 on p.\UThin31 in \cite{BH}.)

This codimensionality result about $I\cdot A[M]$ is new.  It was proved
recently by Gaffney and Massey \cite{GM, (5.7)}, \cite{M98, 4.2} when
$A$ is the ring of germs of a complex-analytic variety, and their proof
inspired ours.  They introduced the remarkable idea of expressing $F$ as
the union of closed sets, each the exceptional divisor of a suitable
blowup, and our corresponding blowup is a stylized version of theirs.
They constructed and used germs of complex-analytic curves in a
remarkable way, which inspired our work with ``paths.''  Their proof and
ours differ mainly because we need to pay careful attention to the
dimension theory, which is so much more delicate for general Noetherian
rings than for geometric rings.  In particular, we must introduce a
certain blowup $B$, dominating $C$, which is unnecessary in their proof.

In the application of the theorem to projective geometry, $X$ is a
variety, $Y$ is contained in its singular locus, and $C$ is its conormal
variety; the latter is the closure of the locus of pairs $(x,H)$ where
$x$ is a simple point and $H$ is a hyperplane tangent to $X$ at $x$.  So
$F$ is a locus of limit tangent hyperplanes at singular points of $X$.
Part~(1) of the theorem provides two cases where $F$ has codimension 1:
Case~(a) $X$ is a singular local complete intersection; Case~(b) the
normal module $\cN$ is not integrally dependent on the Jacobian module
$\cM$.  (In fact, Case~(b) includes Case~(a).)
Thus we obtain a nontrivial lower bound on the dimension of the
dual variety $X'$; see the second author's paper \cite{T}.  Put
differently, when $X'$ is small, the conclusion of (1) fails if $Y$ is
nonempty, and then $\cN$ is dependent on $\cM$.

Part~(2) implies that $F$ is equidimensional if $X$ is a local complete
intersection.  Gaffney and Massey recently proved a similar statement in
complex-analytic geometry.  They applied it to perfect some work of
Gaffney and the first author's in the equisingularity theory of a family
of germs of isolated complete-intersection singularities (ICIS germs),
equipped with a function $f$.  The final result is a definitive
characterization of Thom's A$_f$-condition in terms of the constancy of
numbers of vanishing cycles, or Milnor numbers.

\dsc2 A Global Extension.  It is straightforward, but tedious, to extend
the theorem, obtaining the following corollary, which recovers (10.2)
and (10.3) of \cite{KT}.  A proof will be given in (2.10).

Let $X$ be a Noetherian scheme of finite dimension, $\cN$ a coherent
sheaf, and $\cM$ a proper coherent subsheaf.  Set $Y:=\Supp(\cN/\cM)$.
Let $\cO_X[\cN]$ be a graded quasi-coherent algebra generated by $\cN$
in degree 1, and let $\cO_X[\cM]$ be the subalgebra generated by $\cM$.
Set $P:=\Proj(\cO_X[\cN])$, let $p\:P\to X$ denote the structure map,
and set $Z:=\IV(\cM\cdot \cO_X[\cN])$.  Set $C:=\Proj(\cO_X[\cM])$, let
$c\:C\to X$ denote the structure map, and set $E:=c^{-1}p(Z)$ and
$F:=c^{-1}Y$.  Finally, set $r:=\dim P$.

\pro Corollary Preserve the notation above.
 \part1
 If either \(a)~$\cN$ is locally free of constant rank, $\cO_X[\cN]$ is
the symmetric algebra, $\dim Y<\dim X$, and there exists a point $y\in
Y$ where $\dim\cO_{X,y}=\dim X$, or
 \(b)~$X$ is a closed subscheme of a universally catenary and
biequidimensional scheme, $\dim p^{-1}p(Z)<r$, and
 $Z$ meets an $r$-dimensional component of $P$,
or \(c)~$Z=p^{-1}p(Z)$ as sets, $Z$ is nonempty, $\dim Z<r$, and $X$ is
local, or \(d)~$\dim\cO_{P,z}=r$ for some point $z\in Z$, and $\dim
p^{-1}p(Z)<r$, then
  $$\dim C=r, \and \dim E=r-1.$$
 Furthermore, if $\cN$ is locally free and $\cO_X[\cN]$ is the symmetric
algebra, then $p(Z)=Y$ and $E=F$.
 \part2 Assume either that \(a) holds or that $Z=p^{-1}p(Z)$ as sets,
$\dim Z<r$, and $P$ is equidimensional.  If $X$ is universally catenary
and biequidimensional, then so are $C$, $P$, and $E$.  \endpro

\dsc3 The Height Inequality.  The usual height inequality is this:
        $$d\le m-n+1,\eqno\Cs3.1)$$
 where $d$ is the height of any minimal prime of the ideal of maximal
minors of an $n$ by $m$ matrix with $n\le m$ and with entries in an
arbitrary Noetherian ring $A$.  The inequality is trivial if $d=0$.
Otherwise, as we are now going to see, it results from (1) with (a) of
our theorem (and is nearly equivalent to it); compare \cite{KT, (10.4)}.
Indeed, localizing at the prime and dividing by an arbitrary minimal
prime of $A$, we may assume that $A$ is a local domain of dimension $d$.

Let $M$ be the column space of the matrix, and in the natural way, view
$M$ as a subspace of the free module $N$ of rank $n$.  Then $N/M$ is
supported precisely at the closed point of $\Spec(A)$.  Let $A[N]$ be
the symmetric algebra, $A[M]$ the subalgebra generated by $M$.  Set
$P:=\Proj(A[N])$ and $C:=\Proj(A[M])$.  Standard dimension theory
implies that $\dim P=d+n-1$.  Since $M$ is generated by $m$ elements,
$C$ is a closed subscheme of $\IP_A^{m-1}$.  So $\dim E\le m-1$ since
$E$ is the closed fiber of $C$.  Hence, (1) with (a) of \Cs1) implies
the inequality $d+n-2\le m-1$, and so \Cs3.1).

The height inequality \Cs3.1) can be rewritten in the following form:
        $$m\ge d+n-1.\eqno\Cs3.2)$$
 As such, it is a lower bound on the minimal number $m$ of generators of
a proper submodule $M$ of a free module $N$ over a Noetherian ring,
given in terms of the rank $n$ of $N$ and the height $d$ of any prime
minimal in $\Supp(N/M)$, provided this set is nonempty and nowhere
dense.

The lower bound \Cs3.2) also holds in this general setup: {\it let $A$
be a universally catenary Noetherian ring, $N$ be a finitely generated
module, and $M$ a proper submodule; let $A[N]$ be an arbitrary graded
$A$-algebra generated by $N$ placed in degree $1$, and $A[M]$ the
subalgebra generated by $M$; let $W$ be the subset of $\Spec(A)$ of
primes $\bf p$ where $A[N]_{\bf p}$ is not a finitely generated module
over $A[M]_{\bf p}$; let $\bf p$ be minimal in $W$; and let $\bf q$ be a
homogeneous prime of $A[N]$ such that its contraction ${\bf q}_0:=A\cap
{\bf q}$ is strictly contained in $\bf p$ and the localized quotient
$(A[N]/{\bf q})_{\bf p}$ is not a finitely generated module over
$A[M]_{\bf p}$; then \Cs3.2) holds with, for $m$, the minimal number of
generators of $M$, for $d$, the height of ${\bf p}/{\bf q}_0$, and for
$n$, the transcendence degree of $A[N]/{\bf q}$ over $A/{\bf q}_0$.}
This assertion results similarly from (1) with (b) of \Cs1), after
localizing at ${\bf p}$ and replacing $A[N]$ by $A[N]/{\bf q}$; see the
proof of (1.2) given in (2.10).

\dsc4 Integral Dependence.  Turned around, the corollary yields the
following general criterion for integral dependence in terms of
dimensions.  Let $A$ be a universally catenary and biequidimensional
Noetherian ring, $N$ a finitely generated module, and $M$ a nonzero
proper submodule.  Let $A[N]$ be {\it any\/} graded algebra generated by
$N$ in degree 1, and $A[M]$ the subalgebra generated by $M$.  Define the
maps $p\:P\to X$ and $c\:C\to X$ and the subschemes $Z$ and $E$ as in
\Cs1); set $r:=\dim P$.  Call $N$ {\it integrally dependent\/} on $M$ if
$A[N]$ is a finitely generated module over $A[M]$ (even if $A[N]$ is
{\it not\/} the Rees algebra); it is equivalent to require $Z$ to be
empty, see the middle of (2.1).  Then (1) with (b) of \Cs2) yields this
{\smc Criterion}: {\it $N$ is integrally dependent on $M$ if \(i)~$P$ is
equidimensional, if \(ii)~$\dim p^{-1}p(Z)<r$, and if \(iii)~$\dim E
<r-1$}.

The preceding criterion of ours for modules generalizes the following
criterion of B\"oger's for ideals \cite{B, p.~208}: {\it in a universally
catenary and equidimensional Noetherian local ring $A$, let $M$ and $N$
be nonzero proper ideals with $M\subset N$; then $N$ is integrally
dependent on $M$ if $(\alpha)$~$N_{\bf p}$ is integrally dependent on
$M_{\bf p}$ for every minimal prime $\bf p$ of $A/M$, and
$(\beta)$~$\hgt(M)=\as(M)$ where $\as(M)$ is the analytic spread.}

Indeed, let $A[N]$ and $A[M]$ be the (ordinary) Rees algebras.  Then $C$
and $P$ are the blowups of $\Spec(A)$ along $\IV(M)$ and $\IV(N)$.  So
$C$, $P$, and $X$ are equidimensional of dimension $r$.
Hypothesis~($\alpha$) implies that $p(Z)$ is nowhere dense in $\IV(M)$; so
        $$\dim p(Z) < \dim\IV(M)=r-\hgt(M),$$
 and $\dim p^{-1}p(Z)<r$.  If $\Phi$ denotes the closed fiber of
$C$, then by definition $\as(M):=\dim\Phi+1$.  Hence, standard dimension
theory and  Hypothesis~($\beta$) yield
        $$\dim E\le\as(M)-1+\dim p(Z)<r-1.$$
 Thus all three hypothesis of our criterion hold.

B\"oger replaced Hypothesis~($\alpha$) by the equality of multiplicities,
        $$e(M_{\bf p})=e(N_{\bf p}),$$
 but the two versions of the hypothesis are equivalent by a celebrated
theorem of Rees's.  The latter was generalized to submodules of a free
module by Rees in 4.1 of \cite{R89} and then generalized further
independently by Kirby and Rees in 6.5 of \cite{KR} and by the authors
in (6.7a)(iii) of \cite{KT}.  Also, B\"oger assumed that $A$ is
quasi-unmixed (or formally equidimensional), but this hypothesis implies
that $A$ is universally catenary and equidimensional; see p.\UThin251
and following in \cite{Mat90}.

\dsc5 Projective Geometry. (See \cite{K94} and \cite{T}.) Let $X$ be a
subvariety (or closed, reduced and irreducible subscheme) of dimension
$d$ of the projective $m$-space $\IP^m$ over an algebraically closed
ground field of arbitrary characteristic.  Let $\cI$ be the sheaf of
ideals, and form the usual right exact sequence,
     $$\cI/\cI^2 @>\delta>> \Omega^1_{\IP^m}|X\to \Omega^1_X \to 0.
        \eqno\Cs5.1)$$
 Locally $\delta$ carries a function $f$ vanishing on $X$ to its
differential $df$.  So locally the transpose $\delta^*$ is represented
by a usual Jacobian matrix.

Consider the following nested sequence of three torsion-free sheaves:
        $$\cM:=\Image(\delta^*)\,\,\subset\,\,\cN':=(\Image\delta)^*
                    \,\,\subset\,\,\cN:=(\cI/\cI^2)^*.$$
 where $\cN'$ and $\cN$ are the duals.  The latter is known as the {\it
normal module}.  The first sheaf $\cM$ can be viewed locally as the
column space of a Jacobian matrix; so $\cM$ is known as the {\it
Jacobian module} of $X$.

Let $x\in X$.  First, suppose $X$ is smooth at $x$.  Then \Cs5.1)
splits at $x$.  Hence all three sheaves are equal and are free of rank
$m-d$ at $x$.  Next, suppose $X$ is a complete intersection at $x$.
Then $\cI/\cI^2$ is free at $x$.  So since $X$ is reduced, $\delta$ is
injective.  Hence $\cN'$ and $\cN$ are equal and free at $x$.  If they
are also equal at $x$ to $\cM$, then $x$ must be a simple point because
then $\cI/\cI^2$ is free and \Cs5.1) splits at $x$.  Finally, suppose
$X$ is normal at $x$.  Then $\cN'$ and $\cN$ are equal (but not
necessarily free) at $x$ because $\cO_{X,x}$ satisfies Serre's
conditions (S$_2$) and (R$_1$).

Let $C(X)$ be the {\it conormal variety\/}: by definition, $C(X)$ is the
closure of the set of pairs $(x,H)$ where $x$ is a simple point of $X$
and where $H$ is a hyperplane tangent to $X$ at $x$.  Then $\dim C(X) =
m-1$.  Furthermore,
        $$C(X)=\Proj(\cO_X[\cM])\eqno\Cs5.2)$$
 where $\cO_X[\cM]$ is the Rees algebra (Gaffney, private comm., May
1990).  Indeed, this algebra is sheaf of domains, so
$\Proj(\cO_X[\cM])$ is irreducible.  There is a natural embedding of
the Proj in the product of $\IP^m$ and its dual space: this embedding
is induced by the global Jacobian map,
        $$\cO_X^{m+1}\To\cN(-1),$$
 which arises from the first map $\delta$ in \Cs5.1) and the natural
inclusion map of $\Omega^1_{\IP^m}|X$ into $\cO_X^{m+1}(-1)$.  Finally,
the two sides of \Cs5.2) are equal over the smooth locus of $X$ as
$\cM$ is locally the column space of a Jacobian matrix.

Assume that $x$ is an isolated singular point of $X$ (see \cite{T} for a
more general discussion).  Let $F$ be the fiber of $C(X)$ over $x$.
Part (1) of our theorem in (1.1) implies this: {\it $\dim F=m-2$ either
if $X$ is a complete intersection at $x$ or if the normal module is not
integrally dependent on the Jacobian module at $x$.}

Let $X'$ be the {\it dual variety\/}: by definition, $X'$ is the image
of $C(X)$ under the second projection.  So $X'$ contains the image of
$F$, which may be identified with $F$.  If $X'$ is not a hypersurface
and if $\dim F=m-2$, then $X'=F$.  If the characteristic is zero, then
the dual variety of $X'$ is equal to $X$ (see I-(4) in \cite{K85}).
However, the dual variety $F'$ of $F$ is a cone in $\IP^m$; its vertex
is $x$, and its base is the dual of $F$, viewed as a subvariety of the
hyperplane of hyperplanes through $x$.  Moreover, since $x$ is an
isolated singular point, if $X$ is a cone, then the base is smooth, and
$x$ is the only singular point.  In sum, we have proved this: {\it In
characteristic zero, if the dual variety $X'$ is not a hypersurface and
if, at the isolated singular point $x$, either $X$ is a complete
intersection or, more generally, the normal module is not integrally
dependent on the Jacobian module, then $X'$ has codimension $2$ and $X$
is a cone over a smooth base}.

\dsc6 Example.  The discussion in \Cs5) leads to the following
construction of an example where the conclusion of Part (1) of the
theorem in \Cs1) fails and there is nontrivial integral dependence.
Over an algebraically closed field $k$ of any characteristic, let $G$ be
a smooth subvariety of $\IP^{m-1}$ whose dual variety $G'$ is of
dimension at most $m-3$; specific $G$ will be described below (and more
possible $G$ are described in \cite{K76, p.~360} and \cite{K85, I-7}).
Let $X$ be the projecting cone over $G$ with vertex $x$ in $\IP^m$.
Then its dual variety $X'$ is equal to $G'$.  Hence, by \Cs5), at $x$
the normal module $\cN$ must be integrally dependent on the Jacobian
module $\cM$.  However, algebraically this example is trivial if the two
modules are equal; this possibility will now be investigated.

Sequence \Cs5.1) induces the following short exact sequence:
        $$0\to\cM\to\cN'\to\Ext^1(\Omega^1_X,\cO_X)\to0. $$
 Since $x$ is an isolated singular point, the $\Ext^1$ is concentrated
at $x$.  Moreover, as is well known (see (1.4.3) in \cite{KL} for
example), it is then equal to the module $T^1:=T^1(\cO_X/k,\cO_X)$ of
deformation theory.  Hence, $x$ is a rigid singularity if and only if
$\cM$ and $\cN'$ are equal at $x$.  Moreover, $\cN'$ and $\cN$ are
equal if $G$ is arithmetically normal.

To be specific, let $A$ be an arithmetically normal smooth projective
variety of dimension $a\ge1$, and take $G:=A\x\IP^b$ with $b>a$.  Embed
$G$ via the Segre embedding in $\IP^{m-1}$ say.  Then the dual variety
$G'$ is swept out by the duals of $b$-planes, so has dimension at most
$a+m-2-b$.  Moreover, $G$ is arithmetically normal.  For instance, take
$A:=\IP^a$.  Then $x$ is rigid by a theorem of Thom, Grauert--Kerner,
and Schlessinger; see (2.2.8) in \cite{KL}.  So $\cN=\cM$.  Finally,
take $A$ to be a smooth quartic surface in $\IP^3$, a K3-surface.  Then
the proof of the latter theorem shows that $T^1\neq0$; indeed,
$h^1(\Hom(\Omega^1_A, \cO_A)) =20$ and $h^2(\cO_A)=1$, whence $H^1(\wt
N)\neq0$ where $\wt N$ appears at the end of the proof of (2.2.6) in
\cite{KL}.  So $\cN$ is integrally dependent on $\cM$, but not equal to
it.  This is the desired example.

\dsc7 Equisingularity Theory.  Let $X$ be a complex-analytic germ at $0$
in $\IC^a\x\IC^b$.  Say $X:f_1=0,\dots,f_k=0$  on a neighborhood of $0$ in
$\IC^a\x\IC^b$, where each $f_i$ is an
analytic function $f_i(x,y)$ of the two sets of variables,
        $$x=(x_1,\dots,x_a)\and y=(y_1,\dots,y_b).$$

 For fixed $y$, let $X_y\subset\IC^a$ denote the locus of $x$ such that
$(x,y)\in X$.  Let $Y$ be the locus of $y$ with $(0,y)\in X$, assume
that $Y$ contains a neighborhood of $0$ in $\IC^b$, and identify $Y$
with $0\x Y$.  View $Y$ as the parameter space and $X$ as the total
space of the family of $X_y$.  Finally, assume that the $X_y$ are germs
of isolated complete-intersection singularities (ICIS germs) of
dimension $a-k\ge1$.

Let $f$ be a nonconstant analytic function on $X$ vanishing on $Y$.  Set
   $$\textstyle Z_u:=f^{-1}u \and Z_{u,y}:=Z_u\bigcap X_y.$$
 Let $\Sigma(f)$ denote the ``critical set,'' the union of the singular sets
of the various $Z_u$.  Let $\Sigma_Y(f)$ denote the union of the
singular sets of the various $Z_{u,y}$.

Form the following three conormal varieties: first $C(X,f)$, the closure
of the set of pairs $(w,H)$ where $w:=(x,y)$ is a point of $X-\Sigma(f)$
 and $H$ is a hyperplane in $\IC^a\x\IC^b$ tangent at $w$ to $Z_{fw}$;
second, $C(X,f;Y)$, the closure of the set of $(w,H)$ where $w$ is a
point of $X-Y$ and $H$ is a hyperplane in $\IC^a\x y$ tangent at
$w:=(x,y)$ to $Z_{fw,y}$; third $C(Y)$, the set of pairs
$(w,H)$ where $w:=(0,y)$ is a point of $Y$ and $H$ is a hyperplane
containing $Y$, (in other words, $C(Y)$ is simply $Y\x\IP^{a-1}$).

Extend $f$ over a neighborhood of $X$ in $\IC^a\x\IC^b$ on which
$f_1,\dots,f_k$ are defined, and denote the extension too by $f$; the
choice of extension is immaterial.  Form the following two Jacobian
modules on $X$: first $\cN$, the column space of the Jacobian matrix of
the functions $f_1,\dots,f_k,f$ with respect to all $a+b$ variables
$x,y$; second $\cM$, that with respect to $x$ alone.  So
        $$\cM\subset\cN\subset\cE:=\cO_X^{k+1}.$$

 The reasoning in the proof of \Cs5.2) yields these identifications:
        $$C(X,f)=\Projan(\cN)\and C(X,f;Y)=\Projan(\cM).$$
 Finally, denote the preimage in $C(X,f)$ of $Y$ by $F$, and that of $0$
by $C(X,f)_0$.

Thom's A$_f$-{\it condition\/} at $0$ may be put succinctly as the
condition that
        $$C(X,f)_0\subset C(Y).$$
 It is a well-known preliminary condition for the pair $X,f$ to be
topologically trivially along $Y$ at $0$.  Recently, it was proved to be
equivalent to a weaker condition of topological equisingularity, which
involves the constancy of numbers of vanishing cycles, or Milnor
numbers.

The L\^e-Saito theorem is a celebrated  step in this direction, and
asserts the following: {\it in the case where $X$ is all of\/
$\IC^a\x\IC^b$ and each $Z_{0,y}$ has an isolated singularity at $0$, if
the Milnor number $\mu(Z_{0,y})$ is constant in $y$, then A$_f$ holds.}
L\^e and Saito proved the theorem using Morse theory, but Teissier
reproved it right away using more algebraic-geometric methods.

Following in Teissier's footsteps, Gaffney, Massey and the first author
recently generalized the L\^e-Saito theorem as follows: {\it in the
setup above, at $0$, the germs of $\Sigma(f)$ and $Y$ are equal and
A$_f$ holds if and only if, for $y$ near $0$, the germ $Z_{0,y}$ has an
isolated singularity at $0$, and both $\mu(X_y)$ and $\mu(Z_{0,y})$ are
constant in $y$.} Indeed, \cite{GK, \S5} contains a proof that A$_f$
implies the constancy, and a proof of a weak converse.  The definitive
converse is proved in \cite{GM, (5.8)}; see also \cite{K98, (2.2)}.

A clean composite sketch will now be made of these proofs, highlighting
the use of the complex-analytic version of the theorem in \Cs1).

Suppose that, for $y$ near $0$, the germ $Z_{0,y}$ has an isolated
singularity at $0$, and let $e(y)$ denote the Buchsbaum--Rim
multiplicity of the restriction $\cM|X_y$ in $\cE|X_y$ at $0$.  Theorems
of L\^e and Greuel and of Buchsbaum and Rim yield
        $$e(y)=\mu(X_y)+\mu(Z_{0,y}).$$
 Since Milnor numbers are upper semicontinuous \cite{Lo84, bot.~p.~126},
the two of them are constant in $y$ if and only if $e(y)$ is so.  Thus
we have to prove that, near $0$, the germ $Z_{0,y}$ has an isolated
singularity and $e(y)$ is constant if and only if, at $0$, the germs of
$\Sigma(f)$ and $Y$ are equal and A$_f$ holds.  We'll prove that each
condition holds if and only if, at $0$, the germs of $\Sigma(f)$ and
$Y$ are equal and $\cN$ is integrally dependent on $\cM$.

Suppose A$_f$ holds at $0$.  The gradients of $f_1,\dots,f_k,f$ define
hyperplanes $H$ tangent to the $Z_u$.  So, along any path to $0$ not
lying entirely in $Y$, each $H$ approaches a hyperplane that contains
$Y$.  Therefore, each of the last $b$ components of each gradient
vanishes at $0$ along the curve to order higher than the order of one,
or more, of the first $a$ components.  Hence, by the curve criterion,
$\cN$ is integrally dependent on $\cM$ at $0$.

Suppose that $\Sigma(f)=Y$ and that $\cN$ is integrally dependent on
$\cM$.  Then $\Sigma_Y(f)=Y$; indeed, $\Sigma_Y(f)=\Supp(\cE/\cM)$, and
$\Supp(\cE/\cM)= Y$ as $\cE=\cN$ off $Y$ and a free module is not
dependent on any proper submodule.  Hence $Z_{0,y}$ has an isolated
singularity at $0$.  Let $\cI$ be the zeroth Fitting ideal of $\cE/\cM$.
Then $\cO_X/\cI$ is determinantal, and hence Cohen--Macaulay by a
theorem of Eagon.  Moreover, the support of $\cO_X/\cI$ is equal to $Y$;
in particular, $\cO_X/\cI$ is a finitely generated $\cO_Y$-module.
Consequently, $\cO_X/\cI$ is a Cohen--Macaulay $\cO_Y$-module, and
therefore, by the Auslander--Buchsbaum formula, a free $\cO_Y$-module.
Hence, the restriction of $\cO_X/\cI$ to $X_y$ has constant length.
Since $X_y$ is Cohen--Macaulay, it follows from some theorems of
Buchsbaum and Rim that this length is equal to $e(y)$.  Thus $e(y)$ is
constant.

Conversely, suppose $Z_{0,y}$ has an isolated singularity at $0$.  Then,
replacing $X$ by a smaller representative of its germ, we may assume
that $\Sigma_Y(f)$ is finite over $Y$.  Suppose $e(y)$ is constant.
Then $\Sigma_Y(f)= Y$ because of the upper semi-continuity of the
following sum: the sum, over all the points $w$ in the fiber of
$\Sigma_Y(f)$ over $y$, of the Buchsbaum--Rim multiplicity of the
restriction $\cM|X_y$ in $\cE|X_y$ at $w$.

Hence $\Supp(\cN/\cM)\subset Y$.  So, if $W$ denotes the locus where
$\cN$ is not integrally dependent on $\cM$, then $W\subset Y$.  In fact,
$W\neq Y$ because A$_f$ holds generically on $Y$ by the generic Thom
lemma, and because A$_f$ implies dependence.  Since $\cM$ has $a$
generators, $C(X,f;Y)$ embeds in $X\x\IP^{a-1}$.  Hence, the preimage of
$W$ in $C(X,f;Y)$ has dimension at most $a+b-2$.  However, $C(X,f;Y)$
has dimension $a+b$.  Therefore, by the complex-analytic version of the
criterion of integral dependence discussed in \Cs4) above, $\cN$ is
dependent on $\cM$.

Again suppose that $\Sigma(f)=Y$ and that $\cN$ is integrally dependent
on $\cM$.  Then, as noted above, $\Supp(\cE/\cM)= Y$.  Hence, by the
complex-analytic version of the corollary in \Cs2), each component $F'$
of $F$ has dimension $a+b-1$.  Since $\cN$ is dependent on $\cM$, the
inclusion of $\cM$ into $\cN$ induces a finite surjective map,
  $$g\:C(X,f)\to C(X,f;Y).$$
 Hence, $\dim g(F')=a+b-1$.  However, $C(X,f;Y)\subset X\x\IP^{a-1}$ as
noted above.  Hence $g(F')=Y\x\IP^{a-1}$.  Therefore, $F'$ maps onto
$Y$.  By the generic Thom lemma, the inclusion $F'\subset C(Y)$ holds
generically over $Y$; hence, it holds globally over $Y$.  Thus A$_f$
holds, and the proof is complete.

\sct2. Proof of the theorem and corollary

\dsc1 Preliminaries.  Until the last section \Cs10), preserve the
notation of (1.1), and assume $Y\neq X$.  Form the natural commutative
diagram
        $$\CD B  @>b>> P\\
           @V qVV  @V VpV\\
              C  @>c>> X\endCD$$
 where $B:=\Bl_Z(P)$ is the blowup along $Z$; see \cite{KT, (2.1)}.
Then $Z\neq P$ since $p(Z)\neq X$ and $p$ is surjective.  Hence
$b\:B\to P$ is proper and surjective, as $P$ is integral.  Set
$D:=b^{-1}Z$.

For each nonzero element $\nu$ of $N$, form the ring of elements of
degree 0 and its affine scheme, which is a standard open subscheme of
$P$:
        $$A[N/\nu]:=A[N]_{(\nu)}=A[N]/(1-\nu)\and
          P_\nu:=\Spec(A[N/\nu])\subset P. $$
 Since $\nu$ need not lie in $M$, the corresponding ring and scheme must
be defined differently, but compatibly when $\nu\in M$:
 $$\displaylines{A[M/\nu]:=\Image\big(A[M]\to A[N]\to A[N/\nu]\big);\cr
        Q_\nu:=\Spec(A[M/\nu])\subset Q:=\Spec(A[M]).\cr}$$

The augmentation homomorphism $A[M]\to A$ defines a section of $Q/X$.
Form the corresponding blowups:
        $$\wt Q:=\Bl_XQ\and \wt Q_\nu:=\Bl_{X\cap Q_\nu}Q_\nu.$$
 The exceptional divisor of $\wt Q$ is $C:=\Proj(A[M])$; whence, that of
$\wt Q_\nu$ is
        $$E_\nu:=\wt Q_\nu\cap C.$$
 Since $Q_\nu$ is closed in $Q$, also $E_\nu$ is {\it closed\/} in $C$.

 Note that $E_\nu\subset E:=c^{-1}p(Z)$.  Indeed, work off $p(Z)$, or
assume for the moment that $Z$ is empty.  Then the homogeneous ideal
$M\cdot A[N]$ is irrelevant.  So, if $N_k$ denotes the $k$th graded
piece of $A[N]$, then $M\cdot N_k=N_{k+1}$ for $k\gg0$.  Hence $A[N]$
is a finitely generated $A[M]$-module.  (For use elsewhere, note that
this argument is reversible (compare with \cite{KT, (2.3)}): if $A[N]$
is finitely generated, then $Z$ is empty.)  Hence $A[N/\nu]$ is a
finitely generated $A[M/\nu]$-module.  Now, $M$ generates the unit
ideal in $A[N/\nu]$.  Hence, $M$ generates the unit ideal in
$A[M/\nu]$.  Therefore $X\cap Q_\nu$ is empty off $p(Z)$; whence, so is
$E_\nu$.  Thus $E_\nu\subset E$.

 Set $Z_\nu:=Z\cap P_\nu$.  The inclusion $A[M/\nu]\into A[N/\nu]$
induces maps,
        $$P_\nu\to Q_\nu\and q_\nu\:B_\nu\to\wt Q_\nu
          \where B_\nu:=\Bl_{Z_\nu}(P_\nu)=b^{-1}P_\nu.$$
 Set $D_\nu:=B_\nu\cap D$.  Then the restriction $q_\nu|D_\nu$ is equal
to the restriction,
        $$q\:D_\nu\to E_\nu.$$
 Hence, if $\nu$ varies so that $Z=\bigcup Z_\nu$, then
        $$q(D)=\bigcup q(D_\nu)\subset\bigcup E_\nu\subset E
          \subset F\eqno\Cs1.1)$$

\lem2 If $P$ has dimension $r$ at some $z\in Z_\nu$ for some $\nu$, then
      $$\dim E=\dim F=\dim E_\nu=r-1.$$
 Furthermore, $E_\nu$ is biequidimensional if $A$ is universally
catenary. \endpro

Indeed, the map $P_\nu\to Q_\nu$ carries $z$ to a point $x$ of $X\cap
Q_\nu$, and $x$ must be the unique closed point since $Z$ is closed in
$P$.  Also, $E_\nu$ is the exceptional divisor of $\wt Q_\nu:=\Bl_{X\cap
Q_\nu} Q_\nu$.  Hence, $\dim E_\nu=r-1$ will hold by (3.2)(iii) of
\cite{KT} if
  $$\dim\cO_{Q_\nu,x}=r,\eqno\Cs2.1)$$
 and this equation will now be established.

Generically, the modules $N$ and $M$ are equal to each other, since
$Y\neq X$.  So, generically, the algebras $A[N/\nu]$ and $A[M/\nu]$ are
equal to each other.  Denote their common transcendence degree over $A$
by $f$.  Then $f$ is the dimension of the generic fiber of $p\:P\to X$
as $P_\nu$ is an open subset of $P$.  Set $d:=\dim X$.  Then (3.2)(ii)
of \cite{KT} yields $d+f=r$.

Let ${\bf m}$ be the maximal ideal of $A[M/\nu]$ representing $x$, and
${\bf n}$ that of $A[N/\nu]$ representing $z$.  Then ${\bf n}$ contracts
to ${\bf m}$.  Also, the residue field extension $k({\bf m})/k(x)$ is
trivial.  So, by standard theory \cite{Mat, (14.C), p.~84},
        $$\hgt{\bf m}\leq d+\trdg_AA[M/\nu]-\trdg_{k(x)}k({\bf m})
          =d+f-0=r.$$

By the Hilbert Nullstellensatz, $k({\bf n})/k({\bf m})$ is algebraic.
So, similarly,
  $$\hgt{\bf n}\leq\hgt{\bf m}+\trdg_{A[M/\nu]}A[N/\nu]
    -\trdg_{k({\bf m})}k({\bf n})=\hgt{\bf m}+0-0.$$
 Now, $\hgt{\bf n}=r$ since $P$ has dimension $r$ at $z$.  Hence
$\hgt{\bf m}=r$.  Thus \Cs2.1) holds, and so $\dim E_\nu=r-1$.

Note that $\dim C >\dim F$, for $C$ is irreducible and $C\neq F$ as
$Y\neq X$.  So
        $$r=\dim C >\dim F\ge\dim E\ge\dim E_\nu=r-1$$
 by (1.1.1),  by \Cs1.1), and by what was just proved.  Hence all the
dimensions are as asserted.

Finally, suppose $A$ is universally catenary.  Then $\cO_{Q_\nu,x}$ is
too.  Now, to prove that $E_\nu$ is biequidimensional, we may replace
$Q_\nu$ by $\Spec\cO_{Q_\nu,x}$.  After this replacement, $\wt Q_\nu$
is biequidimensional by (3.8) of \cite{KT}.  Hence its Cartier divisor
$E_\nu$ is too.  The proof is now complete.

\dsc3 Proof of \(1) in\/ {\(1.1)}.  We'll prove that each of the three
hypotheses (a)--(c) implies Hypothesis~(d), that $P$ has dimension $r$
at some $z\in Z$.  Then $z\in Z_\nu$ for some $\nu$ because, as $\nu$
runs through a set of generators of $N$, the various $P_\nu$ cover $P$.
Hence \Cs2) will yield the dimension assertions.

First, assume (a).  Then $A[N]/(M\cdot A[N])$ is equal to the symmetric
algebra on $N/M$, and so $Z=\IP(N/M)$.  Hence $p(Z)=Y$, and so $E=F$.
Hence the closed fiber of $Z$ contains a closed point $z$ because $Y$ is
nonempty as $M\neq N$ by hypothesis.  Finally, $P$ has dimension $r$ at
$z$ by \cite{KT, (3.6)}.

Second, assume $A[N]$ is not a finitely generated module over $A[M]$.
Then $Z$ is nonempty by virtue of part of the argument in \Cs1) showing
$E_\nu\subset E$.  So $p(Z)$ contains the closed point of $X$.  Hence
the closed fiber of $Z$ contains a closed point $z$.  Finally, $P$ has
dimension $r$ at $z$ by \cite{KT, (3.8)}.

Third, assume that $Z=p^{-1}p(Z)$ as sets and that $Z$ is nonempty.
Now, $P$ has dimension $r$ at some point $z$.  Then $p(z)$ is the closed
point of $X$.  So $p(z)\in p(Z)$ since $Z$ is nonempty and is closed.
Hence $z\in Z$ since $Z=p^{-1}p(Z)$.  The proof is now complete.

\dsc4 Paths.  Let $V$ be an $X$-scheme.  By a {\it path\/} to $v\in V$
will be meant an $X$-map $\Spec(R)\to V$, where $R$ is a local
overdomain of $A$, such that the closed point of $\Spec(R)$ maps to $v$.

A path to $w\in C$ is given by a map of graded $A$-algebras $A[M]\to
R[{\bf t}]$ where ${\bf t}$ is an indeterminate; so the path is determined by
the piece in degree 1 of this map, which is an $A$-linear map,
        $$\pi\:M\to R.$$
 Such a $\pi$ will be called a {\it parameterized path}, or {\it pp} for
short.  Of course, $R\pi(M)=R$.

Let $K$ and $L$ be the fraction fields of $A$ and $R$.  Let $A[M]\ox
K\to L[{\bf t}]$ be a map of graded $K$-algebras, $\rho\:M\ox K\to L$
the piece in degree 1.  Suppose $R\rho(M)=Rr$ for some nonzero $r\in L$.
Set ${\bf u}:=r{\bf t}$ and $\pi:=(\rho/r)|M$.  Then $\pi$ is the piece
in degree 1 of the induced map of graded $A$-algebras $A[M]\to R[{\bf
u}]$.  The latter defines a map $\Spec(R)\to C$.  This map is a path to
$w\in C$, where $w$ is the image of closed point of $\Spec(R)$; so $w$
is determined by the composition $M\to R\to k(R)$, where $k(R)$ is the
residue field.

Let $\pi\:M\to R$ be a pp to $w\in C$, and $\pi_K\:M\ox K\to L$ the
extension.  Recall that $M\ox K=N\ox K$.  Suppose $R\pi_K(N)=Rt$ for
some $t\in L$.  Then, by the discussion above, $\psi:=(\pi_K/t)|N$ is a
pp to some $z\in P$.  This pp lifts to a path to some $u\in B$ since
$R\psi(M)=R/t$, and $u\in D$ if $1/t$ lies in the maximal ideal ${\bf
m}_R$.  Furthermore, $q(u)=w$; see \cite{KT, (2.1)}.

For instance, suppose that $A[M]$ is the Rees algebra (or equivalently,
that $A[N]$ is the Rees algebra).  Then $A[M]\ox K$ is the symmetric
algebra over $K$ on the vector space $M\ox K$.  Hence any $A$-linear map
$\rho\:M\ox K\to L$ extends to a map of graded $K$-algebras $A[M]\ox
K\to L[{\bf t}]$.  So, if $R\rho(M)=Rr$, then $\pi:=\rho/r|M$ is a pp to
some $w\in C$, and if $R\pi_K(N)=Rt$ where $t\in {\bf m}_R$, then $w\in
q(D)$.

\lem5 Let $w\in C$.  Then $w\in q(D)$ if there is a pp $\pi\:M\to R$ to
$w$ where $R$ is a valuation ring and if either \(a)~$R\pi_K(N)\neq R$,
or \(b)~$A[N]$ is the Rees algebra, and there is a pp $\theta\:N\to R$ to
a point $z$ of $P$ in $Z$.  \endpro

Indeed, since $R$ is a valuation ring, then $R\pi_K(N)=Rt$ for some
$t\in L$.  If (a) holds, then $1/t\in {\bf m}_R$; whence, $w\in q(D)$
by \Cs4).

Suppose (a) fails, but (b) holds.  By hypothesis, $\theta\:N\to R$ is a
pp to a point of $P$ in $Z$; so $R\theta(N)=R$ and $\theta(M)\subset
{\bf m}_R$.  Since $R$ is a valuation ring, $R\theta(M)=Rs$ for some
$s\in{\bf m}_R$.  If $Rs\neq{\bf m}_R$, then
        $$s=rr'\hbox{ for some }r,r'\in{\bf m}_R;$$
 in fact, any $r\in({\bf m}_R-Rs)$ works.  If $Rs={\bf m}_R$, then the
displayed equation can be achieved by adjoining a square root $r$ of $s$
to $K$ and then replacing $R$ by a valuation ring of the extension
dominating $R$.

The displayed equation implies that $\theta(M)\subset Rs\subset{\bf
m}_Rr$.  Set $\pi':=\pi+(\theta/r)$.  Then $\pi'\:M\to R$ is a pp to
$w\in C$ by \Cs4) because $A[N]$ is the Rees algebra, $R\pi(M)=R$ and
$\theta(M)\subset {\bf m}_Rr$.  Finally, (a) holds for $\pi'$.
Otherwise, (a) would fail for both $\pi$ and $\pi'$.  Then $\theta/r$
would carry $N$ into $R$, and so $\theta(N)\subset Rr\subset {\bf m}_R$.
However, $R\theta(N)=R$ since $\theta$ is a pp.  Replace $\pi$ by
$\pi'$.  Then $w\in q(D)$ by Case~(a).  The proof is now complete.

\lem6 If $\nu$ varies so that $Z=\bigcup Z_\nu$, then
        $$q(D)=\bigcup E_\nu.$$\endpro

Indeed, $q(D)\subset\bigcup E_\nu$ by \Cs1.1).  Conversely, given $w\in
E_\nu$, let $R$ be a valuation ring dominating the local ring of $\wt
Q_\nu$ at $w$, and form the corresponding ring map $\mu\:A[M/\nu]\to R$.
Then $R\mu(M/\nu)=Rr$ for some $r\in{\bf m}_R$ because $E_\nu$ is the
exceptional divisor of $\wt Q_\nu:=\Bl_{X\cap Q_\nu}Q_\nu$.  Moreover,
$\mu$ induces a map of graded $K$-algebras $A[M]\ox K\to L[{\bf t}]$,
whose piece in degree 1 is the composition,
  $$\rho\:M\ox K\to (M/\nu)\ox K\to L.$$
  Then $R\rho(M)=Rr$.  Set $\pi:=\rho/r$.  Then $\pi$ is a pp to $w\in
C$ by \Cs4).  Moreover, $\rho(\nu)=1$, so $R\pi_K(N)\neq R$ as $r\in{\bf
m}_R$.  Hence $w\in q(D)$ by \Cs5) with (a).  The proof is now complete.

\dsc7 About $q(D)$. If one of the hypotheses (a)--(d) of (1) in (1.1) holds,
then \Cs6) and the proof in \Cs3) of (1) imply that
        $$\dim q(D)=r-1.$$
  Similarly, if $A$ is universally catenary, then $q(D)$ is
biequidimensional, being the union of closed biequidimensional subsets,
one for each $\nu$ such that the closed fiber of $Z_\nu$ is nonempty.

\lem8 If \(a)~$N$ is free and $A[N]$ is the symmetric algebra, or
\(b)~$M$ is free and $A[M]$ is the symmetric algebra, or
\(c)~$Z=p^{-1}p(Z)$ as sets, then $q(D)=E$.\endpro

Indeed, suppose (c) holds.  Then $D=q^{-1}E$ as sets.  Now, $q\:B\to C$
is surjective by \cite{KT, (2.6)} since $A[M]\subset A[N]$.  Hence
$q(D)=E$.

Suppose (a) or (b) holds.  Then $C\x_XP$ is integral, and it dominates
both $C$ and $P$.  Let $w\in E$.  Since $E:=c^{-1}p(Z)$, then
$c(w)=p(z)$ for some $z\in Z$.  Let $v\in C\x_XP$ map to both $w$ and
$z$, and let $R$ be a valuation ring dominating the local ring at $v$.
Then the natural map $\Spec(R)\to C$ is a path to $w$, and the natural
map $\Spec(R)\to P$ is a path to $z$.  Hence $w\in q(D)$ by \Cs5) with
(b).  Thus, $q(D)\supset E$, and the converse inclusion holds by
\Cs1.1).  The proof is now complete.

\dsc9 Proof of \(2) in\/ {\(1.1)}.  Assume $A$ is universally catenary.
Then $C$ and $P$ are biequidimensional by \cite{KT, (3.8)} since they
are irreducible and the maps $p\:P\to X$ and $c\:C\to X$ are proper.
Furthermore, $q(D)$ is biequidimensional by \Cs7).  Finally, $q(D)=E$ by
\Cs8) if $N$ is free or if $Z=p^{-1}p(Z)$ as sets.  The proof is now
complete.

\dsc10 Proof of the corollary in\/ {\(1.2)}.  Preserve the notation of
(1.2), and form the natural commutative diagram
        $$\CD B  @>b>> P\\
           @V qVV  @V VpV\\
              C  @>c>> X\endCD$$
 where $B:=\Bl_Z(P)$ is the blowup along $Z$; see \cite{KT, (2.1)}.
Then $q\:B\to C$ is surjective by \cite{KT, (2.6)} since
$\cO_X[\cM]\subset\cO_X[\cN]$.

Suppose for a moment that $\cN$ is locally free of rank $n$ and that
$\cO_X[\cN]$ is the symmetric algebra.  Then $Z=\IP(\cN/\cM)$.  Hence
$p(Z)=Y$, and so $E=F$.  Also, if $\dim Y<\dim X$, then, by \cite{KT,
(3.6)} applied at each closed point of $Y$ and at each of $X$,
        $$\dim p^{-1}Y=\dim Y+n-1<\dim X+n-1=\dim P =:r.$$
 Hence $\dim p^{-1}p(Z)<r$ under Hypothesis~(a) too.

Let $C'$ be a component of $C$.  Since $q$ is surjective, $C'=q(B')$
for some component $B'$ of $B$.  Then, by \cite{KT, (3.2)},
        $$\dim C'\le\dim B'=\dim b(B')\le\dim P=:r.\eqno\Cs2.1)$$
 Let $E'$ be a component of $E$.  Then $\dim E'\le r-1$.  Otherwise,
$E'$ is a component $C'$ of $C$ of dimension $r$ because of \Cs2.1).
Then $cq(B')\subset p(Z)$, so $b(B')\subset p^{-1}p(Z)$.  So, if $\dim
p^{-1}p(Z)<r$, then there is a contradiction.

Consider Part~(1).  To prove $\dim C=r$ and $\dim E=r-1$, it is enough,
by the above, to find one $C'$ of dimension $r$ and one $E'$ of
dimension $r-1$.

We may assume (d) holds, as there is a point $z\in Z$ with
$\dim\cO_{P,z} =r$
 also if (a), (b) or (c) holds.
  Indeed, if (a) holds,
then any closed point $z\in Z$ lying over $y\in Y$ will do by \cite{KT,
(3.6)} applied locally at $y\in X$.  If (b) holds, then $Z$ meets an
$r$-dimensional component $P'$ of $P$, and any closed point $z\in P'\cap
Z$ will do by \cite{KT, (3.8)}.  Finally, assume (c), and take $z\in P$
such that $\dim\cO_{P,z} =r$.  Then $p(z)$ is the unique closed point of
$X$.  So $p(z)\in p(Z)$ since $Z$ is nonempty and is closed.  Hence
$z\in Z$ since $Z=p^{-1}p(Z)$.

We may localize the setup at $p(z)$.  Then $X$ is the spectrum of a
local ring, $A$ say.  Moreover, the two $\cO_{X}$-algebras are
associated to two graded $A$-algebras, $A[N]$ and $A[M]$ say.

Take a component $P'$ of $P$ whose local ring at $z$ has dimension $r$,
and give $P'$ its reduced structure.  Then $P'=\Proj(A'[N'])$ where
$A'[N']$ is a graded domain and a quotient of $A[N]$.  Let $M'$ be the
image of $M$ in $N'$, and use a prime to indicate the corresponding
constructions.  Then $Z'=Z\cap P'$; hence (d) continues to hold.
Moreover, the $C'$ and $E'$ are closed subsets of  $C$ and $E$; hence,
the latter have  dimensions $r$ and $r-1$ if the former do.

By hypothesis, $\dim p^{-1}p(Z)<r$; so $p^{-1}p(Z')\neq P'$ and so
$p(Z')\neq X'$.  Since $A'$ is a domain and $N'$ is torsion free, $N'$ is
free on a dense open set, $U$ say, of $X'$.  Then $Y'\cap U=p(Z')\cap U$
by the second paragraph of the proof.  Hence $Y'\neq X'$.  So
$\dim C'=r$ by (1.1.1), and $\dim E'=r-1$ by (1) with (d) of the
theorem.  The proof of (1) of the corollary is now complete.

Consider Part~(2).  To prove it, we may localize at an arbitrary closed
point of $X$.  Indeed, since $X$ is universally catenary and
biequidimensional, $P$ is equidimensional also if (a) holds by (3.2)(ii)
of \cite{KT}, and so $P$ is biequidimensional by (3.8) of \cite{KT}.
We may also replace $P$ by an arbitrary component $P'$.  Indeed, each
component $C'$ of $C$ corresponds to some $P'$ by the third paragraph
above.  Let $E'$ be a component of $E$.  Then $E'$ lies in some $C'$,
which corresponds to some $P'$.  If $p^{-1}p(Z)=Z$, then any closed
point of $P'$ lies in $Z$.  If $(a)$ holds, then the whole closed fiber
of $P$ lies in $P'$, and so any closed point of $Z$ lies in $P'$.
Finally, proceeding as in the proof of Part~(1), reduce (2) of the
corollary to (2) of the theorem.  Thus (2) is proved, and the proof of
the corollary is complete.

\sctn References

\references

B
 E. B\"oger,
 Verallgemeinerung eines Multiplizit\"atensatzes von D. Rees,
 \ja 12 1969 207--15.

BH
 W. Bruns and J. Herzog,
 ``Cohen-Macaulay rings,''
 Cambr. studies in adv. math {\bf 39}, Cambr. U. Press, 1993.

BR
 D. A. Buchsbaum and D. S. Rim,
 A generalized Koszul complex. II. Depth and multiplicity,
 \tams 111 1963 197--224.

BV
  W. Bruns and U. Vetter,
 ``Determinantal Rings,''
 Spr.L.N.M. {\bf 1327}, 1988.

GM
 T. Gaffney and D. Massey,
 Trends in equisingularity theory,
 {\it in} the proceedings of the 1996 Liverpool conference in honor of
CTC Wall, singularities volume, W. Bruce and D. Mond (eds.), Cambr. U. Press.

GK
 T. Gaffney and S. L. Kleiman,
 Specialization of integral dependence for modules,
 {\it alg-geom/}9610003.

KR
 D. Kirby and D. Rees,
   Multiplicities in graded rings I: The general theory,
   {\it in} ``Commutative algebra: Syzygies, multiplicities, and
  birational algebra, Proceedings, Mount
  Holyoke College (USA) 1992," W. J. Heinzer,
  C. L. Huneke, J. D. Sally (eds.),
   \conm 159 1994 209--67.

K76
 S. L. Kleiman,
 The enumerative theory of singularities,
 {\it in} ``Real and complex singularities, Proceedings, Oslo 1976,''
 P.  Holm (ed.), Sijthoff \& Noordhoff 1977, pp. 297--396.

K85
 S. L. Kleiman,
 Tangency and duality,
 {\it in} ``Proc.~1984 Vancouver Conf. in Algebraic Geometry,''
 J.  Carrell, A.  V.  Geramita, P. Russell (eds.),
 CMS Conf.~Proc.~{\bf 6}, AMS,  1986, pp. 163--226.

K94
 S. L. Kleiman,
 A generalized Teissier--Pl\"ucker formula,
 {\it in} ``Classification of algebraic varieties," Ciliberto,
L. Livorni, and A. Sommese (eds.), \conm 162 1994 249--60.

K98
S. L. Kleiman,
 Equisingularity, multiplicity, and dependence,
 {\it in} the volume of proceedings of the June 1997 Ferrara conference
in honor of M. Fiorentini, to be published by Marcel Dekker.

KL
 S. L. Kleiman and J. Landolfi,
 Geometry and deformations of special Schubert varieties,
 \comp 23 1971 407--434.

KT
 S. Kleiman and A. Thorup,
 A geometric theory of the Buchsbaum--Rim multiplicity,
 \ja 167 1994 168--231.

Lo84
 E. J. N.  Looijenga,
 Isolated singular points on complete intersections,
 London Mathematical Society lecture note series {\bf 77},
Cambridge University Press, 1984.

M98
 D. B. Massey,
 Critical points of functions on singular spaces, {\it Preprint\/} 1998.

Mat
 H. Matsumura,
 ``Commutative algebra,''
 W.A. Benjamin, inc., New York, 1970.

Mat90
 H. Matsumura,
 ``Commutative ring theory,''
 Cambr. studies in adv. math {\bf 8}, Cambr. U. Press, 1990.

R89
 D. Rees,
  ``Gaffney's problem,''
  manuscript dated February 22, 1989.

T
 A. Thorup,
 The Teissier--Pl\"ucker formula,
 work in progress.

\endreferences
\bye